\documentclass[11pt,a4paper]{article}
\usepackage[T1]{fontenc}
\usepackage[utf8]{inputenc}
\usepackage{lmodern}
\usepackage{microtype}
\usepackage{amsmath,amssymb,amsthm,mathtools}
\usepackage{geometry}
\usepackage{xcolor}
\usepackage{enumitem}
\usepackage{hyperref}
\usepackage{mathrsfs}
\usepackage{booktabs,tabularx,pdflscape,array}
\hypersetup{colorlinks=true,linkcolor=blue!50!black,citecolor=blue!50!black,urlcolor=blue!55!black}
\allowdisplaybreaks
\numberwithin{equation}{section}
\newtheorem{theorem}{Theorem}[section]
\newtheorem{proposition}[theorem]{Proposition}
\newtheorem{lemma}[theorem]{Lemma}
\newtheorem{corollary}[theorem]{Corollary}
\newtheorem{definition}[theorem]{Definition}
\theoremstyle{remark}
\newtheorem{remark}[theorem]{Remark}

\newcommand{\dd}{\,\mathrm d}
\newcommand{\cV}{\mathcal V}

\title{\bfseries Weak and Sharp Stretched-Exponential Flocking in the Kinetic Cucker-Smale Model with Noncompactly Supported Initial Data}
\author{Ioannis Markou}
\date{}

\begin{document}
	\maketitle
	\begin{abstract}
   We study the long-time behavior of the kinetic Cucker–Smale model under initial data with unbounded support in both position and velocity. Such data invalidate classical diameter-based flocking arguments and make sustained spatial cohesion difficult to establish. We prove an estimate for spatial characteristics under a finite spatial second moment and a sufficiently high velocity moment, without any requirement  on coupling strength. We also prove weak flocking, i.e. control of spatial variance and integrability of velocity fluctuations in time. Under stronger assumptions on the initial tails, we establish uniform propagation of spatial exponential moments of explicitly determined order, a substantially stronger conclusion than bounded spatial variance. These estimates lead to explicit stretched-exponential rates of velocity alignment and convergence of the centered spatial distribution toward a limiting profile. Another contribution is the construction of countably atomic solutions that show optimality of the decay exponent in the spatial-tail-limited class. Together, the results clarify how initial phase-space tails govern spatial confinement and quantitative flocking for kinetic distributions without compact support.
	\end{abstract}
	
	
	\section{Introduction}
	\subsection{From classical flocking to noncompact kinetic distributions}
	Collective alignment is a basic mechanism in the mathematical description of
	self-organized motion. In the Cucker--Smale (CS) model, individuals adjust their
	velocities by averaging the velocities of the surrounding population, with an
	influence that decreases as the spatial separation increases. The original
	particle system of Cucker and Smale \cite{CuckerSmale2007} prompted a large
	literature on unconditional and conditional flocking, mean-field limits, and
	the associated kinetic and hydrodynamic equations. The passage from particle
	dynamics to kinetic descriptions was developed, in particular, by Ha and
	Tadmor \cite{HaTadmor2008}, while the Lyapunov approach of Ha and Liu
	\cite{HaLiu2009} established a useful connection between alignment
	mechanisms and mean-field stability. Kinetic flocking and measure-valued
	well-posedness were further investigated by Carrillo, Fornasier, Rosado,
	and Toscani \cite{CarrilloEtAl2010} and by Ca\~nizo, Carrillo, and Rosado
	\cite{CanizoCarrilloRosado2011}. An exhaustive presentation of the literature on this subject is far beyond the scope of this work, so we give the following surveys \cite{ChoiHaLi2017, Ha2025Survey, MinakowskiMuchaPeszekZatorska2019}, and also the work of the author on the subject \cite{HaskovecMarkou2020a, HaskovecMarkou2020b, Markou2018}.
	
	In this article, we work with the kinetic CS equation 
	\begin{equation}\label{eq:intro-model}
		\partial_t f_t+v\cdot\nabla_x f_t+\nabla_v\cdot(F[f_t]f_t)=0,
		\qquad
		F[f_t](t,x,v)=\kappa\int\psi(|x-y|)(w-v)\,f_t(\dd y,\dd w),
	\end{equation}
	with coupling strength $\kappa>0$, and communication weight
	\begin{equation}\label{eq:intro-kernel}
		\psi(r)=(1+r^2)^{-\beta/2},\qquad 0<\beta<1.
	\end{equation}
	The unknown, $f_t(x,v)$ is a probability measure that describes the probability of  particles located at the phase-space point $(x,v) \in \mathbb R_x^d\times\mathbb R_v^d$ at time $t$.  
	The measure $f_t$ is normalized to have unit total mass. Conservation of
	momentum fixes the mean velocity $v_c$, and the center of mass of the system moves as
	$x_c(t)=x_c(0)+t v_c$. The two natural macroscopic fluctuations are
	\begin{equation}\label{eq:intro-moments}
		\mathcal{X}(t)=\int|x- x_c(t)|^2\,\dd f_t(x,v),\qquad
		\mathcal V(t)=\int|v- v_c|^2\,\dd f_t(x,v).
	\end{equation}
	
	For compactly supported initial data, the traditional flocking argument
	controls the diameter of the spatial support and thereby recovers a strictly
	positive lower bound on the communication weight along interacting pairs.
	This ultimately yields exponential decay of the velocity diameter and a
	uniformly bounded spatial diameter in the unconditional regime. This mechanism of control is not available for fully noncompact data: even when nearly all of the mass stays
	within a small region, particles of arbitrarily small mass can occur at
	arbitrarily large distances. In this setting, the infimum of the pairwise
	communication coefficients may be zero at every time.
	
	This obstruction motivates a distinction between several notions of
	flocking. Velocity alignment in the mean-square sense means
	$\mathcal V(t)\to 0$. The usual \emph{weak flocking} conclusion supplements
	this by $\sup_{t\ge0} \mathcal{X}(t)<\infty$. The latter is an averaged form of
	spatial cohesion but does not, by itself, propagate a prescribed class of
	spatial tails. In particular, boundedness of $\mathcal{X}(t)$ yields only a
	Chebyshev estimate of the type $f_t\{|x- x_c(t)|>R\} \le C R^{-2}$, for any $R>0$; it does not
	imply a uniform exponential or stretched-exponential spatial bound. A
	central objective here is the stronger property
	\begin{equation}\label{eq:intro-uniform-tail}
		\sup_{t\ge0}\int_{\mathbb R^{2d}}
		\exp\big(a_*|x-x_c(t)|^\sigma\big)\,\dd f_t(x,v)<\infty
		\qquad\text{for some }a_*,\sigma>0.
	\end{equation}
	Besides quantifying spatial cohesion, such a bound provides an
	exponentially small error when the velocity dissipation is localized to
	an expanding ball. This connection is essential to the convergence rates
	proved below.
	
	\subsection{Previous results for noncompactly supported data}
	The study of the kinetic model with noncompact support is not new.
	 From the point of view of well-posedness,
	Chen and Yin \cite{ChenYin2019} proved global weak existence under finite
	first phase-space moments. They also showed stronger classical and measure
	well-posedness results for compactly supported data. Their results
	illustrate why one must separate the existence of a noncompact weak
	solution from the propagation of strong flocking estimates. Wang and
	Xue \cite{WangXue2023} studied formation and velocity consensus for
	kinetic CS systems with noncompact spatial distributions, advancing the
	asymptotic theory beyond the usual bounded spatial support setting.
	
	A particularly important sequence of advances is due to S.-Y. Ha, X. Wang et al. In particular, Ha, Wang, and Xue \cite{HaWangXue2025} considered kinetic
	distributions with noncompact spatial support and bounded initial velocity
	support. Under sufficiently rapid spatial decay and an appropriate 
	coupling hypothesis, they obtained weak flocking, namely decay of the
	centered velocity variance together with boundedness of the centered
	spatial second moment. Their analysis uses an infinite-particle
	approximations and tail estimates in place of the classical
	support-diameter functional.
	
 	Ha and Wang \cite{HaWang2026} subsequently treated the \emph{fully
	noncompact} setting, allowing both position and velocity support to be
	unbounded for any $d\geq 1$. They distinguished polynomially
	and exponentially decaying initial distributions and developed a
	\emph{time-varying effective-region} method. Their polynomial-moment
	theorem gives algebraic decay of the velocity variance and uniform
	spatial second-moment control; the moment exponents are linked by
	explicit inequalities. For data with an exponential moment of the form
	$\int e^{a(|x|+|v|)}\,\dd f_{\rm in}<\infty$, they established weak
	flocking with convergence faster than any prescribed algebraic rate.
	If the velocity-tail hypothesis is strengthened to
	$\int e^{a(|x|+|v|^\delta)}\,\dd f_{\rm in}<\infty$, with $\delta>1$
	and $\beta<1-1/\delta$, their Theorem~2.2 gives the explicit estimate
	\begin{equation}\label{eq:intro-Ha-rate}
		\mathcal V(t)\le C\exp\!\left[-ct^{\,1-\beta\delta/(\delta-1)}\right].
	\end{equation}
	The same paper proves weak spatial cohesion. Its stated theorems do not,
	however, assert a uniform-in-time spatial exponential moment of the
	form \eqref{eq:intro-uniform-tail}. This distinction matters: an
	exponential bound on the initial distribution and a uniform second
	spatial moment at later times do not automatically furnish such a
	uniform exponential-moment estimate.
	
	The question of constructing the corresponding solutions is equally
	important. The noncompact weak-existence theory of Chen--Yin
	\cite{ChenYin2019} gives one point of reference, while Ha and Wang
	\cite{HaWang2026} establish uniqueness under additional bounded-density
	and exponential velocity-tail assumptions. More recently, Ha and Wang
	\cite{HaWangStability2026} obtained finite-time Osgood-type stability
	and random mean-field approximation for a class with exponential
	velocity moments and finite spatial second moments. An essential
	technical point is that boundedness and smoothness of $\psi$ do
	\emph{not} give a globally uniform phase-space Lipschitz constant for
	$F[f]$ when $|v|$ is unbounded: the spatial derivative of the
	alignment force involves the factor $w-v$. We therefore prove existence
	of Lagrangian probability solutions in the moment class used below
	and state our estimates for every such solution, without asserting a
	general nonlinear uniqueness theorem.
	
	Other recent contributions clarify adjacent, but distinct, mechanisms.
	Ha and Wang's relativistic model \cite{HaWangRelativistic2025} treats
	noncompact phase-space data with a modified velocity structure.
	Ha, Wang, and Yoon \cite{HaWangYoon2026} investigate difficulties in
	uniform-in-time stability and propose a modified noncompact interaction
	model. Ha and Wang \cite{HaWangConfining2026} establish quantitative
	flocking in the presence of a confining force; genuine exponential
	fluctuation-energy relaxation there is assisted by confinement and
	should not be confused with the free alignment dynamics in
	\eqref{eq:intro-model}. Finally, Ha, Wang, and Xue
	\cite{HaWangXueInfinite2026} analyze clustering and spectral properties
	of a countably infinite particle system. That perspective is relevant
	to the construction of infinitely many distant ``satellite'' particles,
	but their infinite-graph spectral conclusions address a different
	question from the propagation of kinetic spatial exponential moments.
	
	\subsection{Motivation of this work and main results}
	The developments we mentioned leave room for a more detailed description of the
	\emph{spatial distribution itself}. Rather than merely asking whether
	$\mathcal{X}(t)$ remains finite, we ask whether the initial position and velocity
	tails control the location of each characteristic for \emph{all} time,
	and whether that confinement yields sharp quantitative velocity decay.
	The answers presented in this paper rest on three related results.
	
	Our first result is a finite-moment confinement principle, valid for
	every $d\ge1$. Suppose the centered initial spatial second moment is
	finite and the initial centered velocity has a moment of order
	\begin{equation}\label{eq:intro-moment-threshold}
		m>\frac{2}{1-\beta}.
	\end{equation}
	We prove weak spatial cohesion, quantitative algebraic velocity
	alignment, and the integrability of the velocity standard deviation,
	\begin{equation}\label{eq:intro-J}
		\sup_{t\ge0} \mathcal{X}(t)<\infty,
		\qquad \int_0^\infty\sqrt{\mathcal V(t)}\,\dd t<\infty.
	\end{equation}
	The key additional consequence is the \emph{support-independent}
	characteristic estimate
	\begin{equation}\label{eq:intro-confinement}
		\sup_{t\ge0}|X_t(z)|\le
		C\Big(1+|X_0(z)|+|V_0(z)|^{1/(1-\beta)}\Big)
		\quad\text{for }f_{\rm in}\text{-a.e. }z,
	\end{equation}
	where $X_t,V_t$ are centered characteristic variables. The constant
	is controlled by averaged initial moments and fixed model parameters;
	it is independent of the support diameter, particle number, and the
	individual label $z$. While related characteristic damping ideas are
	classical, the combination with the finite-moment Lyapunov estimate
	allows their use for fully noncompact measures in any dimension.
	
	The second result is uniform propagation of strong spatial tails. If
	for some $p,q,a,b>0$,
	\begin{equation}\label{eq:intro-exponential-data}
		\int\exp\big(a|x- x_c(0)|^p+b|v- v_c|^q\big)
		\,\dd f_{\rm in}(x,v)<\infty,
	\end{equation}
	then \eqref{eq:intro-uniform-tail} holds with the explicitly determined
	order
	\begin{equation}\label{eq:intro-sigma}
		\sigma=\min\{p,(1-\beta)q\}.
	\end{equation}
	Thus the spatial tail is limited by the initial position exponent
	$p$ and by the velocity exponent $q$ after multiplication by the
	coercivity exponent $1-\beta$. The statement includes fully noncompact
	phase-space distributions, and neither compact velocity support nor a
	large-coupling condition is assumed. The associated velocity-variance
	estimate is
	\begin{equation}\label{eq:intro-present-rate}
		\mathcal V(t)\le C\exp\left[-c(1+t)^{\sigma/(\sigma+\beta)}\right],
		\qquad t\ge0.
	\end{equation}
	In addition, the centered spatial marginal converges to a limiting
	probability measure in $W_2$, with a quantitative rate. This conclusion
	identifies an asymptotic spatial profile even though the original
	spatial support need never become bounded.
	
	Our third result concerns optimality. For every $p>0$, we construct a
	symmetric countably atomic kinetic solution, supported on a fixed line
	in $\mathbb R^d$, with bounded initial velocities and uniformly bounded
	spatial exponential moments of order $p$, for which
	\begin{equation}\label{eq:intro-sharpness}
		\mathcal V(t)\ge c\exp\left[-C(1+t)^{p/(p+\beta)}\right].
	\end{equation}
	The interaction of an arbitrarily distant satellite with the central
	mass may be extremely weak, even when the satellite has exponentially
	small weight. Choosing a geometric sequence of radii and suitable
	satellite masses realizes this slow mode simultaneously with strong
	spatial-tail bounds. The lower bound shows that the exponent
	$p/(p+\beta)$ cannot be improved uniformly in the spatial-tail-limited
	measure class with bounded velocities. The construction does not claim
	sharpness for smooth strictly positive kinetic densities or for the
	regime in which the velocity-tail exponent determines $\sigma$.
	
	For comparison with the exponential part of Ha--Wang
	\cite[Theorem~2.2]{HaWang2026}, take $\beta=1/4$ and assume initial
	exponential position tails of order one and Gaussian-type velocity
	tails of order two. Their stated estimate yields exponent $1/2$,
	whereas \eqref{eq:intro-sigma}--\eqref{eq:intro-present-rate} yield
	$\sigma=1$ and exponent $4/5$ in arbitrary dimension. Under merely
	exponential velocity tails of order one, the present result instead
	gives $\sigma=3/4$ and exponent $3/4$; the corresponding statement in
	\cite{HaWang2026} asserts super-polynomial alignment but does not
	specify this exponent. These are comparisons between \emph{stated}
	sufficient estimates, not claims that the earlier methods could not
	be sharpened. Our finite-moment theorem also improves the sufficient moment
	conditions for weak flocking. It requires only a finite initial
	spatial second moment and a finite velocity moment of order
	$m>2/(1-\beta)$, for every $0<\beta<1$.
	The polynomial-moment criterion in
	\cite[Theorem~2.2]{HaWang2026} also admits a spatial moment of
	order two in part of the parameter range. More precisely,
	setting their spatial and velocity moment orders to $D_1=2$
	and $D_2=m$, respectively, and optimizing over their auxiliary
	parameters gives the sufficient conditions
	\[
	0<\beta<\frac12,
	\qquad
	m>\frac{2(1-\beta)}{1-2\beta}.
	\]
	For $0<\beta<1/2$, our velocity-moment threshold is strictly
	smaller; for $1/2\leq\beta<1$, their stated polynomial-moment
	criterion cannot be satisfied with $D_1=2$, whereas our theorem
	still applies under the stated velocity-moment assumption.
	Thus the improvement concerns the admissible combinations of
	spatial and velocity moments, rather than the admission of
	finite-second-moment spatial data alone.
	
	The results of the present article consequently contribute in three
	complementary directions: a finite-moment, dimension-independent
	criterion for spatial cohesion and characteristic confinement;
	uniform propagation of quantitative spatial exponential tails with a
	corresponding explicit velocity decay rate; and a countably atomic
	example matching that rate in the spatial-tail-limited class. The
	first result identifies a tractable threshold within our method, but
	we do not claim that the strict moment inequality
	\eqref{eq:intro-moment-threshold} is necessary in every dimension.
	The atomic construction proves sharpness in a specified measure class,
	without settling smooth-density sharpness. These distinctions are
	important for separating proven results from possible extensions.
	
	The rest of the paper is arranged as follows. Section~2 states the
	main results and notation. Section~3 establishes existence of global
	Lagrangian solutions under the required finite moments. Section~4
	records the velocity dissipation identities and localization tools.
	Section~5 proves the coupled Lyapunov estimate and the
	characteristic confinement theorem. Section~6 establishes strong
	tail propagation, stretched-exponential alignment, and convergence
	of the centered spatial marginal. Section~7 develops the countably
	atomic sharpness construction. The concluding section identifies
	remaining mathematical questions.
	
	\section{Setting and main theorems}\label{sec:setting}
	Let $d\ge1$, $0<\beta<1$, and $\kappa>0$. The kinetic CS model (for $(x,v)\in \mathbb{R}^d \times \mathbb{R}^d$) is
	\begin{align} 	\label{eq:model} \begin{cases}
		& \partial_t f_t+v\cdot\nabla_x f_t+\nabla_v\cdot(F[f_t]f_t)=0, \\
		& F[f_t](t,x,v)=\kappa\int
		\psi(|x-y|)(w-v)\,f_t(\dd y,\dd w),
		\qquad \psi(r)=(1+r^2)^{-\beta/2}.
		 \end{cases}
	\end{align}
	Let $f_{\rm in}(x,v)$ be a Borel probability measure of total mass one, i.e $\int \dd f_{\rm in}(x,v)=1$. Conservation of mass implies that any solution remains a probability measure $f_t \in \mathcal{P}(\mathbb{R}_{x,v}^{2d})$. The equation admits solutions that can be represented as the push-forward of the initial measure $f_{\rm in}$ by the transport flow. A global \textit{Lagrangian} solution is a family $f_t=(\mathbf X_t,\mathbf V_t)_\#f_{\rm in}$ with almost-everywhere
	locally absolutely continuous trajectories satisfying the kinetic
	characteristic equations and locally continuous in
	$L^2(f_{\rm in};\mathbb R^{2d})$. The push-forward identity gives, for every bounded Borel function $\phi:\mathbb R^{2d}\to \mathbb R$, or more generally every $\phi\in L^1(f_t)$ that 
	\begin{align*} \int \phi(x,v) \dd f_t(x,v) = \int \phi(\mathbf{X}_t(x,v),\mathbf{V}_t(x,v)) \dd f_{\rm in}(x,v). \end{align*}
	Theorem~\ref{thm:existence} below
	establishes the existence of at least one such solution for every
	initial probability measure with the stated finite moments; all our
	a priori estimates hold for every Lagrangian solution in this class.
	We do not impose or assert general nonlinear uniqueness.
	Lemma~\ref{lem:finite-time} supplies the finite-time identities.
	
	\begin{lemma}[Spatial and velocity averages] We define the following averages for positions and velocities 
		\begin{align*} x_c(t):=\int x \dd f_t , \qquad v_c(t):=\int v \dd f_t .\end{align*}
		Then we have that $v_c(t)=v_c(0)=:v_c$, and $x_c(t)=x_c(0)+v_c t$.
	\end{lemma}
	\begin{proof} The proof follows by simply differentiating $x_c(t)$ and $v_c(t)$ along solutions of \eqref{eq:model} 
	\end{proof}
	
	Before we go to the main theorems we give some definitions for later use, starting from the centered characteristics and position-velocity fluctuations
	\begin{definition}[Centered characteristics and variances] For any point $z=(x,v) \in \mathbb{R}^d \times  \mathbb{R}^d$ we define the centered characteristics from point $z$ by
		\begin{equation}
		X_t(z):=\mathbf X_t(z)- x_c(t),\quad
		V_t(z):=\mathbf V_t(z)-v_c.
		\end{equation}
		We also define the variances around mean position and velocity by
		\begin{equation}\label{eq:SV}
			\mathcal{X}(t):=\int|X_t|^2\,\dd f_{\rm in},
			\qquad \cV(t):=\int|V_t|^2\,\dd f_{\rm in},
			\qquad \mathcal{Z}(t):=\sqrt{1+\mathcal{X}(t)}.
		\end{equation}
		The quantity $\mathcal{Z}(t)$ is a renormalized spatial variance that comes handy later.
	\end{definition}
	It is not difficult to verify that
	\begin{align*}
		\int X_t\,\dd f_{\rm in}&=\int V_t\,\dd f_{\rm in}=0 \end{align*}
		We also have the characteristic system around the average
		\begin{align} \label{eq:characteristic} \dot X_t=V_t, \qquad
		\dot V_t(z)&=-\kappa\int\psi(|X_t(z)-X_t(z')|)
		(V_t(z)-V_t(z'))\,\dd f_{\rm in}(z').
		\end{align}
		We also give the standard definition of the $W_2$ Wasserstein distance which is a metric distance for distributions.
		We denote by $\mathcal P_2(\mathbb R^d)$ the space
		of Borel probability measures with finite second
		moments. For $\rho,\nu\in\mathcal P_2(\mathbb R^d)$,
		the quadratic Wasserstein distance is defined by
		\[
		W_2(\rho,\nu)
		=
		\left(
		\inf_{\pi\in\Pi(\rho,\nu)}
		\int_{\mathbb R^d\times\mathbb R^d}
		|x-y|^2\,d\pi(x,y)
		\right)^{1/2},
		\]
		where $\Pi(\rho,\nu)$ is the set of all couplings
		of $\rho$ and $\nu$.

	In our first theorem we prove weak-flocking under the assumptions of quadratic spatial initial moments, and velocity moments of order $m>2/(1-\beta)$. We also show an important global (in time) estimate of the spatial characteristics.   
	\begin{theorem}[Weak flocking and control of spatial characteristic]
		\label{thm:finite-moment}
		Assume $\mathcal{X}(0)<\infty$ and for some $m>2/(1-\beta)$, we assume finite $m$'th velocity moment of the initial data, i.e.
		\begin{equation}\label{eq:finite-m-assumption}
			M_m:=\int|v-v_c|^m\,\dd f_{\rm in}<\infty.
		\end{equation}
		Then the following estimates hold:
		\begin{align}
			\sup_{t\ge0}\mathcal{X}(t)&<\infty,\label{eq:finite-S}\\
			\int_0^\infty\sqrt{\cV(t)}\,\dd t&<\infty,\label{eq:int-sqrtV}\\
			\sup_{t\ge0}|X_t(z)|&\le C\left(1+|X_0(z)|+
			|V_0(z)|^{1/(1-\beta)}\right)
			\quad\text{for} \quad f_{\rm in}\text{-a.e. }z.\label{eq:finite-pointwise}
		\end{align}
		Moreover, if $\cV(0)>0$, then for
		$\eta=\beta m/[2(m-2)]\in(0,1/2)$,
		\begin{equation}\label{eq:finite-polynomial-rate}
			\cV(t)\le C(1+t)^{-1/\eta}.
		\end{equation}
		The constant in \eqref{eq:finite-pointwise} depends only on
		$d,\beta,\kappa,\mathcal{X}(0),\cV(0),m,M_m$.
		If $\cV(0)=0$, all centered velocities vanish and the conclusions
		are immediate.
	\end{theorem}
	A by-product of Theorem \ref{thm:finite-moment} is that we can actually get $L^2$ and $W_2$ control of the centered spatial tails
	\begin{corollary}[Polynomial spatial convergence under finite moments]
		\label{cor:polynomial-spatial}
		Under the assumptions of Theorem~\ref{thm:finite-moment},
		there exists $X_\infty\in L^2(f_{\rm in};\mathbb R^d)$ such that
		\[
		\|X_t-X_\infty\|_{L^2(f_{\rm in})}+
		W_2((X_t)_\#f_{\rm in},(X_\infty)_\#f_{\rm in})
		\le C(1+t)^{-(\frac1{2\eta}-1)},
		\qquad \eta=\frac{\beta m}{2(m-2)}<\frac12.
		\]
		For $\mathcal V(0)=0$, the centered positions are constant.
	\end{corollary}
	
	\begin{theorem}[Exponential spatial tails and stretched-exponential flocking]
		\label{thm:main}
		Assume, for $p,q>0$ and $a,b>0$,
		\begin{equation}\label{eq:initial-tail}
			M_{a,b,p,q}:=\int
			\exp\left(a|x-x_c(0)|^p+b|v-v_c|^q\right)\,\dd f_{\rm in}<\infty.
		\end{equation}
		Set
		\begin{equation}\label{eq:sigma}
			\sigma=\min\{p,(1-\beta)q\},
			\qquad \vartheta=\frac{\sigma}{\sigma+\beta}.
		\end{equation}
		There exist $a_*,c,C>0$ depending on the initial moments and
		parameters such that
		\begin{align}
			\sup_{t\ge0}\int e^{a_*|x-x_c(t)|^\sigma}\,\dd f_t &\le C,
			\label{eq:main-tail}\\
			\cV(t)&\le C e^{-c(1+t)^\vartheta}.
			\label{eq:main-rate}
		\end{align}
		All assertions hold for arbitrary $d\ge1$ and $\kappa>0$.
	\end{theorem}
	
	\begin{corollary}[Spatial tails, strong convergence and $W_2$]
		\label{cor:W2}
		Under Theorem~\ref{thm:main}, for every finite $s\ge0$,
		\begin{equation}\label{eq:strong-moment-tails}
			\sup_{t\ge0}\int|X_t|^s\,\dd f_{\rm in}<\infty,
			\qquad \lim_{R\to\infty}\sup_{t\ge0}
			\int_{\{|X_t|>R\}}|X_t|^s\,\dd f_{\rm in}=0.
		\end{equation}
		There is $X_\infty\in L^2(f_{\rm in};\mathbb R^d)$ satisfying
		\begin{equation}\label{eq:W2-rate}
			\|X_t-X_\infty\|_{L^2(f_{\rm in})}+
			W_2((X_t)_\#f_{\rm in},(X_\infty)_\#f_{\rm in})
			\le C e^{-c(1+t)^\vartheta},
		\end{equation}
		and the limiting spatial measure has an exponential moment of order
		$\sigma$ (for the same $a_*$ after weakening it if necessary).
	\end{corollary}
	
	\begin{remark} We mention the fact that in both corollaries above it suffice to prove the inequalities for $	\|X_t-X_\infty\|_{L^2(f_{\rm in})}$ due to the inclusion
		\begin{align*} 	W_2((X_t)_\#f_{\rm in},(X_\infty)_\#f_{\rm in})^2 \leq \int |X_t(z)-X_{\infty}(z)|^2\, \dd f_{\rm in}(z) . \end{align*}
	The inequality follows from the definition of $W_2$, and the fact that $\pi_t=(X_t,X_{\infty})_\#f_{\rm in}$ is only one of the possible couplings.
	\end{remark}
	
	\begin{proposition}[Conditional quantitative flocking for $d\geq 1$]
		\label{prop:conditional}
		Suppose $d\ge1$ and $\beta>0$. If a global centered Lagrangian
		solution satisfies
		\begin{equation}\label{eq:conditional-assumptions}
			\sup_{t\ge0}\int e^{a|X_t|^\sigma}\,\dd f_{\rm in}\le M_x<\infty,
			\qquad \int|V_0|^{2+\varepsilon}\,\dd f_{\rm in}<\infty
		\end{equation}
		for some $a,\sigma,\varepsilon>0$, then
		\begin{equation}\label{eq:conditional-conclusion}
			\cV(t)\le C\exp[-c(1+t)^{\sigma/(\sigma+\beta)}].
		\end{equation}
		This proposition does not use the restriction $0<\beta<1$;
		that restriction enters through Theorem~\ref{thm:finite-moment}.
	\end{proposition}
	
	\section{Existence theory and Lagrangian compatibility}\label{sec:existence}
	\begin{theorem}[Existence for the noncompact initial data]
		\label{thm:existence}
		Let $d\ge1$, $\kappa>0$, and let $\psi$ be bounded, nonnegative,
		globally Lipschitz and s.t. $\|\psi\|_{L^{\infty}} \leq 1$. For every Borel probability measure
		$f_{\rm in}\in\mathcal P(\mathbb R^{2d})$ with
		\[
		\int (|x|^2+|v|^m)\,\dd f_{\rm in}<\infty,\qquad m\ge2,
		\]
		there exists a global weak kinetic solution
		$f\in C([0,\infty);\mathcal P_1(\mathbb R^{2d}))$
		which has a deterministic Lagrangian representation
		$f_t=(\mathbf X_t,\mathbf V_t)_\#f_{\rm in}$. The characteristic map is
		locally continuous in $L^2(f_{\rm in})$, and
		\[
		\sup_{t\le T}\int (|x|^2+|v|^m)\,\dd f_t<\infty
		\qquad (T<\infty).
		\]
		The centered velocity $m$-moment is nonincreasing and the exact
		velocity-variance identity holds. These conclusions hold for every
		such Lagrangian solution, not only for the one constructed below.
	\end{theorem}
	\begin{proof}
		\emph{Step 1: finite atomic approximation.}
		Approximate $f_{\rm in}$ by finite atomic probabilities $f_{\rm in}^N$ obtained by
		truncation and quantization so that $f_{\rm in}^N\to f_{\rm in}$ in $W_2$, their
		initial spatial second moments converge, and their velocity
		$m$-moments converge (truncation can be chosen separately in the two
		variables). Every finite particle system is global. Its velocities
		remain in their initial convex hull, and symmetric pairing of the
		interaction terms gives conservation of mass and momentum, exact
		second-moment dissipation, and nonincrease of each convex centered
		velocity moment. In particular, with constants independent of $N$,
		\begin{equation}\label{eq:existence-uniform}
			\sup_{N,t\le T}\int|v|^m\,\dd f_t^N<\infty,
			\qquad
			\sup_{N,t\le T}\int|x|^2\,\dd f_t^N<\infty
			\quad (T<\infty).
		\end{equation}
		The second estimate follows from
		$\sqrt{\int|x|^2\,\dd f_t^N}\le
		\sqrt{\int|x|^2\,\dd f_{\rm in}^N}+t\sqrt{\int|v|^2\,\dd f_{\rm in}^N}$.
		
		\emph{Step 2: a global weak limit.}
		The phase-space velocities of the particle curves are
		$(v,F[f^N])$, and
		$|F[f^N](x,v)|\le\kappa(|v|+\int|w|\,\dd f^N)$.
		Thus $t\mapsto f_t^N$ is equi-Lipschitz in $W_1$ on each finite time
		interval. Estimate \eqref{eq:existence-uniform} gives tightness and
		uniform integrability of first moments; metric Arzel\`a--Ascoli and a
		diagonal argument yield a subsequence converging in
		$C([0,T];\mathcal P_1)$ for every $T$.
		In the weak kinetic formulation, the nonlinear term is
		\[
		\kappa\iint \nabla_v\varphi(x,v)\cdot
		\psi(|x-y|)(w-v)\,\dd f_t^N(x,v)\dd f_t^N(y,w).
		\]
		The bounded continuity of the kernel and the uniform $m$-moment
		bound ($m\ge2$) give convergence after a velocity truncation; the
		remaining tails are uniformly small because the integrand grows at
		most linearly in $|v|+|w|$. Hence the limiting curve solves the
		kinetic equation globally, has the initial measure $f_{\rm in}$, and
		satisfies the finite-time moment bounds by lower semicontinuity.
		
		\emph{Step 3: deterministic characteristics.}
		For this fixed solution, write
		$F[f](t,x,v)=-\kappa a_t(x)v+\kappa b_t(x)$ with
		\[
		a_t(x)=\int\psi(|x-y|)\,\dd f_t(y,w),\qquad
		b_t(x)=\int\psi(|x-y|)w\,\dd f_t(y,w).
		\]
		If $L=\operatorname{Lip}\psi$ and
		$M_1(t)=\int|w|\,\dd f_t$, then
		\[
		|a_t(x)-a_t(y)|\le L|x-y|,\quad
		|b_t(x)-b_t(y)|\le LM_1(t)|x-y|.
		\]
		Consequently, on a region $|v|\le K$, the full field
		$(v,F[f](t,x,v))$ is locally Lipschitz in $(x,v)$, with a constant
		depending on $K$ and the finite-time first velocity moment.
		\emph{It is not globally Lipschitz in phase space when $v$ is
			unbounded.} The field has at most linear growth and is integrable
		against $f_t$. The superposition principle for the continuity
		equation (Theorem 8.2.1 in \cite{AmbrosioGigliSavare2008})
		 represents $f_t$ as a
		probability on its absolutely continuous characteristic curves.
		Local Lipschitz continuity and linear growth give uniqueness and
		global continuation of the ordinary differential equation for this
		\emph{fixed} $f$, so disintegration by the initial label is a Dirac
		mass. Thus $f_t=(\mathbf X_t, \mathbf V_t)_\#f_{\rm in}$ for a deterministic
		flow; this does not imply uniqueness among distinct nonlinear
		solutions.
		
		\emph{Step 4: regularity and identities.}
		For any such characteristic solution the damping variation-of-
		constants formula first gives, on $[0,T]$,
		$|V_t(z)|\le |V_0(z)|+C_T$. The finite-time moment bound then
		justifies differentiation of the centered velocity variance,
		recovering its exact dissipation identity, and improves this to
		$|V_t(z)|\le |V_0(z)|+\kappa T\sqrt{\mathcal V(0)}$.
		Integration of $\dot X_t=V_t$ yields an $L^2(f_{\rm in})$-dominating
		function for positions and velocities, so the flow is locally
		continuous in $L^2(f_{\rm in})$. The same domination and convex
		truncations establish nonincrease of centered velocity $m$-moments.
		The argument applies to every Lagrangian solution with the stated
		finite-time moments.
	\end{proof}
	\begin{remark}[Relation to earlier solution theory]
		Chen--Yin \cite{ChenYin2019} supply global noncompact weak existence
		under a first-moment hypothesis. Ha--Wang \cite{HaWang2026} and
		\cite{HaWangStability2026} give uniqueness or quantitative weak
		stability under additional density and exponential-velocity
		conditions. Neither result should be read as a blanket uniqueness
		theorem for all the atomic data admitted here. The argument above
		provides existence in the exact Lagrangian class needed for the
		a priori estimates, without such a uniqueness claim.
	\end{remark}
	
	\section{Finite time control of characteristics and main velocity localization estimate}
	\label{sec:prelim}
	We begin with a simple technical lemma that is helpful when we differentiate in time across integrals
	\begin{lemma}[Finite-time control of characteristics and velocity variance identity]
		\label{lem:finite-time}
		Assume $\mathcal{X}(0)<\infty$ and $\cV(0)<\infty$. Then on each finite interval
		$[0,T]$ the characteristic solution satisfies
		\begin{align}
			\sup_{0\le t\le T}|V_t(z)|
			&\le |V_0(z)|+\kappa T\sqrt{\cV(0)},\label{eq:FTU}\\
			\sup_{0\le t\le T}|X_t(z)|
			&\le |X_0(z)|+T|V_0(z)|+
			\frac{\kappa T^2}{2}\sqrt{\cV(0)}.\label{eq:FTX}
		\end{align}
		The velocity variance is locally absolutely continuous, and
		\begin{equation}\label{eq:variance-diss}
		\dot \cV(t)=-\kappa\iint
			\psi(|X_t-X_t'|)|V_t-V_t'|^2\,\dd f_{\rm in}\dd f_{\rm in}'\le 0
			\quad\text{for a.e. }t.
		\end{equation}
		When an initial $m$'th velocity moment is finite, the same is true on
		every compact time interval.
	\end{lemma}
	\begin{proof}
		Write the characteristic equation in damping form
		$\dot V_t=-\kappa a_t(X_t)V_t+\kappa b_t(X_t)$, where
		$0\le a_t\le1$ and $|b_t(x)|\le m_1(t):=\int|V_t|\,\dd f_{\rm in}$.
		The scalar integrating-factor formula, followed by integration in
		$f_{\rm in}$, yields
		\[
		|V_t(z)|\le |V_0(z)|+\kappa\int_0^t m_1(s)\,\dd s,
		\qquad m_1(t)\le m_1(0)e^{\kappa t}.
		\]
		The stipulated local $L^2$ continuity ensures that the quantities
		in this initial argument are locally finite. For $t\le T$ these
		inequalities give a pointwise bound $|V_t(z)|\le |V_0(z)|+C_T$;
		hence the second moment, and any initially finite higher polynomial
		velocity moment, admit a uniform integrable dominator on $[0,T]$.
		We may consequently differentiate $\int|V_t|^2\,\dd f_{\rm in}$
		and symmetrize, obtaining \eqref{eq:variance-diss}. Thus
		$\cV(t)\le\cV(0)$, whence
		$|b_t(x)|\le m_1(t)\le\sqrt{\cV(0)}$.
		Returning to the integrating-factor formula gives
		\eqref{eq:FTU}; integrating $\dot X_t=V_t$ gives
		\eqref{eq:FTX}. For a higher convex polynomial test function,
		its differentiated integrand is dominated on $[0,T]$ by a constant
		times $1+|V_0(z)|^m+|V_0(z')|^m$, which permits direct dominated
		convergence, or equivalently a truncation argument.
	\end{proof}
	 With the help of the above lemma we get
	\begin{lemma}[Convex velocity moments]
		\label{lem:convex}
		For a convex differentiable $G:\mathbb R^d\to\mathbb R$ of
		suitable polynomial growth, the map
		$t\mapsto\int G(V_t)\,\dd f_{\rm in}$ is nonincreasing. In particular,
		for $m\ge2$,
		\begin{equation}\label{eq:propagate-m}
			\int|V_t|^m\,\dd f_{\rm in}\le\int|V_0|^m\,\dd f_{\rm in}.
		\end{equation}
	\end{lemma}
	\begin{proof}
		The identity
		\begin{align*}
			\frac{\dd}{\dd t}\int G(V_t)\,\dd f_{\rm in}
			=-\frac\kappa2\iint\psi(|X_t-X_t'|)
			(\nabla G(V_t)-\nabla G(V_t'))\cdot(V_t-V_t')\,
			\dd f_{\rm in}\dd f_{\rm in}'
		\end{align*}
		is nonpositive by monotonicity of the gradient of a convex
		function. The previous finite-time bounds justify the identity
		for $G(V)=|V|^m$ by truncation.
	\end{proof}
	
	\begin{lemma}[Localized velocity coercivity]\label{lem:localized}
		For every $R>0$, we define
		\[
		B_R(t):=\{|X_t|\le R\},\quad Q_R(t):=f_{\rm in}(B_R(t)^c),\quad
		T_R(t):=\int_{B_R(t)^c}|V_t|^2\,\dd f_{\rm in}.
		\]
		Then for almost every $t$,
		\begin{equation}\label{eq:localized}
		\dot \cV(t)\le -2\kappa\psi(2R)
			\left[(1-Q_R(t))\cV(t)-T_R(t)\right].
		\end{equation}
	\end{lemma}
	\begin{proof}
		On $B_R \times B_R$, the communication weight is at least
		$\psi(2R)$. The exact identity
		\[
		\iint_{B_R\times B_R}|V_t-V'_t|^2\,\dd f_{\rm in}\dd f_{\rm in}'
		=2(1-Q_R)\int_{B_R}|V_t|^2\,\dd f_{\rm in}
		-2\left|\int_{B_R}V_t \,\dd f_{\rm in}\right|^2
		\]
		and the zero-mean condition give
		$|\int_{B_R}V_t \,\dd f_{\rm in}|^2
		=|\int_{B_R^c}V_t \,\dd f_{\rm in}|^2\le Q_RT_R$.
		Consequently the restricted double integral is bounded below by
		$2[(1-Q_R)\cV-T_R]$. Combining this with
		\eqref{eq:variance-diss} proves \eqref{eq:localized}.
	\end{proof}
	
	\section{Weak flocking and spatial characteristic confinement}
	\label{sec:multidimensional}
	In the result that follows, we put to use the coercivity estimate we gave in Lemma~ \ref{lem:localized} to prove an inequality for the dissipation of the velocity variance in terms of the spatial and velocity variances. Then we introduce a Lyapunov functional to prove control of spatial variance, and show weak flocking.
	\begin{lemma}[The Lyapunov functional]\label{lem:coupled}
		Assume the hypotheses of Theorem~\ref{thm:finite-moment}, and
		$\cV(0)>0$. Set
		\begin{equation}\label{eq:eta}
			\gamma=\frac{m}{2(m-2)},\qquad
			\eta=\beta\gamma=\frac{\beta m}{2(m-2)},\qquad
			\delta=\frac12-\eta>0.
		\end{equation}
		There exists $k>0$ depending only on the initial moments and
		$\kappa,\beta,m$ such that
		\begin{equation}\label{eq:nonlinear-diss}
		\dot \cV(t)\le-k\mathcal{Z}(t)^{-\beta}\cV(t)^{1+\eta}
		\end{equation}
		for almost every $t$ with $\cV(t)>0$. Moreover, the functional
		\begin{equation}\label{eq:global-Lyapunov}
			\mathscr L(t):=\mathcal{Z}(t)^{1-\beta}+
			\frac{1-\beta}{k\delta}\cV(t)^\delta
		\end{equation}
		is dissipative, i.e.  $\dot{\mathscr L}(t)\leq 0$.
	\end{lemma}
	
	\begin{proof}
		The spatial moment gives
		\begin{equation}\label{eq:Zprime}
			\dot{\mathcal{Z}}(t)=\frac{\dot{\mathcal{X}}(t)}{2 \mathcal{Z}(t)}=\frac{\int X_t\cdot V_t\,\dd f_{\rm in}}{\mathcal{Z}(t)}
			\le\sqrt{\cV(t)}.
		\end{equation}
		At a time where $\cV(t)>0$, choose
		\begin{equation}\label{eq:adaptive-radius}
			R(t)=A \mathcal{Z}(t)\left(\frac{\cV(0)}{\cV(t)}\right)^\gamma,
		\end{equation}
		with $A\ge2$ fixed below. Since $\cV(t)\le\cV(0)$,
		Chebyshev's inequality and \eqref{eq:propagate-m} give
		\begin{align}
			Q_{R(t)}(t)
			&\le A^{-2}\left(\frac{\cV(t)}{\cV(0)}\right)^{2\gamma},
			\label{eq:Q-adaptive}\\
			T_{R(t)}(t)
			&\le M_m^{2/m}Q_{R(t)}(t)^{1-2/m}
			\le M_m^{2/m}A^{-2(m-2)/m}
			\frac{\cV(t)}{\cV(0)}.\label{eq:T-adaptive}
		\end{align}
		Choose $A$ so large that $A^{-2}\le1/4$ and
		$M_m^{2/m}A^{-2(m-2)/m}/\cV(0)\le1/4$.
		Then $Q_{R(t)}\le1/4$ and $T_{R(t)}\le\cV(t)/4$.
		Since \eqref{eq:localized} follows pointwise from the variance
		identity at each of its differentiability times and is valid for
		\emph{every} $R>0$ at such a time, substituting the time-dependent
		radius causes no moving-boundary differentiation issue.
		Thus \eqref{eq:localized} becomes
		\begin{equation}\label{eq:pre-nonlinear}
		\dot \cV(t)\le-\kappa\psi(2R(t))\cV(t).
		\end{equation}
		Writing $\Lambda=(\cV(0)/\cV(t))^\gamma\ge1$ and using
		$\mathcal{Z}\ge 1$,
		\begin{align*}
			\psi(2R(t))
			&=(1+4A^2 \mathcal{Z}^2\Lambda^2)^{-\beta/2}\\
			&\ge(1+4A^2)^{-\beta/2}\mathcal{Z}^{-\beta}
			\left(\frac{\cV(t)}{\cV(0)}\right)^{\eta}.
		\end{align*}
		Thus \eqref{eq:nonlinear-diss} holds with
		\begin{equation}\label{eq:k-constant}
			k=\kappa(1+4A^2)^{-\beta/2}\cV(0)^{-\eta}>0.
		\end{equation}
		Differentiating the two terms in \eqref{eq:global-Lyapunov},
		and using \eqref{eq:Zprime} and
		$\delta+\eta=1/2$, gives
		\begin{align*}
			\dot{\mathscr L}(t)
			\le(1-\beta) \mathcal{Z}(t)^{-\beta}\sqrt{\cV(t)}
			-\frac{1-\beta}{k\delta}\,
			k\delta \mathcal{Z}(t)^{-\beta}\cV(t)^{\delta+\eta}=0.
		\end{align*}
		The estimates extend through possible zero-variance times by the
		absorbing character of perfect alignment, or by limiting arguments.
	\end{proof}
	
	\begin{corollary}[Spatial cohesion and integrable velocity variance]
		\label{cor:poly}
		Under the hypotheses of Theorem~\ref{thm:finite-moment},
		$\mathcal{Z}(t)\le \mathcal{Z}_*$ uniformly in time. In particular,
		\begin{equation}\label{eq:polynomialV}
			\dot \cV(t)\le-k_*\cV(t)^{1+\eta},
			\qquad \cV(t)\le
			\bigl(\cV(0)^{-\eta}+\eta k_*t\bigr)^{-1/\eta},
		\end{equation}
		where $k_*=k \mathcal{Z}_*^{-\beta}>0$. Moreover,
		\begin{equation}\label{eq:sqrtV-integrable}
			\int_0^\infty\sqrt{\cV(t)}\,\dd t<\infty.
		\end{equation}
	\end{corollary}
	\begin{proof}
		Monotonicity of \eqref{eq:global-Lyapunov} gives the explicit bound
		\[
		\mathcal{Z}(t)\le \mathcal{Z}_*:=\left(\mathcal{Z}(0)^{1-\beta}+
		\frac{1-\beta}{k\delta}\cV(0)^\delta\right)^{1/(1-\beta)}.
		\]
		Integrating \eqref{eq:nonlinear-diss} with $\mathcal{Z}(t)\le \mathcal{Z}_*$
		proves \eqref{eq:polynomialV}. Since $\eta<1/2$,
		$\cV(t)^{1/2}\lesssim(1+t)^{-1/(2\eta)}$ is integrable;
		indeed the polynomial estimate gives the explicit bound
		\[
		\int_0^\infty\sqrt{\cV(t)}\,\dd t
		\le \frac{\cV(0)^\delta}{k_*\delta}<\infty,
		\qquad \delta=\frac12-\eta>0.
		\]
	\end{proof}
	
	\begin{proof}[Proof of Corollary~\ref{cor:polynomial-spatial}]
		For $t\le s$, the fundamental theorem of calculus and Minkowski's
		inequality yield
		$\|X_s-X_t\|_{L^2(f_{\rm in})}\le\int_t^s\sqrt{\mathcal V(\tau)}\,\dd\tau$.
		The bound \eqref{eq:polynomialV} implies
		\[
		\int_t^\infty\sqrt{\mathcal V(\tau)}\,\dd\tau
		\le C(1+t)^{1-1/(2\eta)}.
		\]
		The Cauchy criterion gives $X_t\to X_\infty$ strongly in
		$L^2(f_{\rm in})$ with that rate, and the canonical coupling gives the
		same bound in $W_2$.
	\end{proof}
	
	We now proceed to the proof of the spatial characteristic estimate
	\begin{lemma}[Uniform control of spatial characteristic]
		\label{lem:pointwise}
		Under the hypotheses of Theorem~\ref{thm:finite-moment}, we get
		\begin{equation}\label{eq:pointwise-confinement}
			\sup_{t\ge0}|X_t(z)|\le
			C\left(1+|X_0(z)|+|V_0(z)|^{1/(1-\beta)}\right)
		\end{equation}
		for $f_{\rm in}$-a.e. $z$, with $C$ controlled by the initial averaged
		moments and $\kappa,\beta$.
	\end{lemma}
	\begin{proof}
		By \eqref{eq:finite-S}, choose $R_0<\infty$ such that
		$f_{\rm in}\{|X_t|\le R_0\}\ge1/2$ for almost every $t>0$. The characteristic velocity equation is
		\begin{equation}\label{eq:AB}
			\dot V_t=-\kappa a_t(X_t)V_t+\kappa b_t(X_t), \end{equation} for
			\begin{equation*}
			\quad a_t(x)=\int\psi(|x-X_t'|)\,\dd f_{\rm in}',
			\quad b_t(x)=\int\psi(|x-X_t'|)V_t'\,\dd f_{\rm in}'.
		\end{equation*}
		For every $x\in\mathbb R^d$ and $t\ge0$,
		\begin{equation}\label{eq:A-lower}
			a_t(x)\ge\tfrac12\psi(|x|+R_0)
			\ge c_0(1+|x|)^{-\beta},\qquad
			c_0=\tfrac12(1+R_0)^{-\beta},
		\end{equation}
		and
		\begin{equation}\label{eq:B-upper}
			|b_t(x)|\le\int|V_t'|\,\dd f_{\rm in}'
			\le\sqrt{\cV(t)}.
		\end{equation}
		Fix a label $z$ and put $r(t)=|X_t(z)|$,
		$u(t)=|V_t(z)|$. The upper Dini derivatives satisfy
		\begin{equation}\label{eq:ru-ineq}
			D^+u(t)\le-\kappa c_0(1+r)^{-\beta}u+
			\kappa\sqrt{\cV(t)},
			\qquad D^+r(t)\le u(t).
		\end{equation}
		The use of upper Dini derivatives defined by 
		\begin{align*} 	D^+u(t):=\limsup \limits_{h \downarrow 0}\frac{u(t+h)-u(t)}{h}\end{align*}
		accounts for the possible vanishing of $r(t)$ or $u(t)$. At nonzero velocity, the first estimate follows by taking the scalar product of the velocity equation with the unit velocity vector. At zero velocity, its upper right Dini derivative is bounded by just the magnitude of the forcing term. Similarly, the position equation implies that the upper right Dini derivative of the spatial radius does not exceed the velocity magnitude.
		
		Now, for $\alpha=1-\beta$,
		\[
		D^+(1+r)^\alpha\le\alpha(1+r)^{-\beta}u.
		\]
		Therefore the damping terms cancel in the combination
		\begin{equation}\label{eq:characteristic-Lyapunov}
			D^+\left[u(t)+\frac{\kappa c_0}{\alpha}(1+r(t))^\alpha\right]
			\le\kappa\sqrt{\cV(t)}.
		\end{equation}
		Integrating, by \eqref{eq:sqrtV-integrable}, gives
		\begin{equation}\label{eq:char-conf-explicit}
			(1+r(t))^\alpha\le(1+r(0))^\alpha+
			\frac{\alpha}{\kappa c_0}u(0)+
			\frac{\alpha}{c_0}\int_0^\infty\sqrt{\cV(s)}\,\dd s.
		\end{equation}
		The additivity and elementary power inequalities yield
		\eqref{eq:pointwise-confinement}. The constant $c_0$ is
		independent of the characteristic label; thus the final coefficient
		is independent of any support diameter.
	\end{proof}
	
	\begin{proof}[Proof of Theorem~\ref{thm:finite-moment}]
		If $\cV(0)=0$, then $V_t=0$ and $X_t=X_0$ almost everywhere;
		the conclusions follow. Otherwise Lemma~\ref{lem:coupled}
		and Corollary~\ref{cor:poly} give the spatial second-moment and
		velocity-integrability conclusions, and Lemma~\ref{lem:pointwise}
		gives characteristic confinement. Equation
		\eqref{eq:finite-polynomial-rate} follows from
		\eqref{eq:polynomialV}.
	\end{proof}
	
	\section{Exponential tails and stretched-exponential velocity variance}
	\label{sec:tails}
	Having established the uniform spatial characteristic control inequality, we can use it to prove propagation of exponential tails 
	\begin{proof}[Proof of the spatial-tail assertion of Theorem~\ref{thm:main}]
		The initial exponential velocity moment provides a finite moment
		of order $m>2/(1-\beta)$, while the initial exponential spatial
		moment provides $\mathcal{X}(0)<\infty$. Thus
		Theorem~\ref{thm:finite-moment} applies. It gives
		\[
		\sup_{t\ge0}|X_t(z)|\le
		C\bigl(1+|X_0(z)|+|V_0(z)|^{1/(1-\beta)}\bigr).
		\]
		For $\sigma=\min\{p,(1-\beta)q\}$, there is $C_1$ such that
		\begin{equation}\label{eq:exp-domination}
			|X_t(z)|^\sigma\le
			C_1(1+|X_0(z)|^p+|V_0(z)|^q).
		\end{equation}
		Choosing $a_*>0$ with $a_*C_1\le\min\{a,b\}$, taking the exponent 
		and integrating proves \eqref{eq:main-tail}. 
	\end{proof}
	
	Next comes the proof of stretched-exponential velocity variance. First we prove it in Proposition~\ref{prop:conditional}, and then in Theorem~\ref{thm:main}
	\begin{proof}[Proof of Proposition~\ref{prop:conditional}]
	 Assume that for some $\sigma,a>0$ and $M_x<\infty$,
		\begin{equation}\label{eq:spatial-moment}
			\sup_{t\ge0}\int e^{a|X_t(z)|^\sigma}\dd f_{\rm in}(z)\le M_x,
		\end{equation}
		and that for some $r=2+\varepsilon>2$,
		\begin{equation}\label{eq:vel-moment}
			M_r:=\int|V_0(z)|^r\dd f_{\rm in}(z)<\infty.
		\end{equation}
		Then there exist constants $c,C>0$ such that
		\begin{equation}\label{eq:target}
			\cV(t) \le C\exp\!\left[-c(1+t)^{\sigma/(\sigma+\beta)}\right],\qquad t\ge0.
		\end{equation}
	
		\textbf{Step 1: The spatial mass and velocity energy outside a ball.}
		For any $R>0$ and $t\ge0$, remember that
		\[
		B_R(t)=\{z:|X_t(z)|\le R\},\qquad
		Q_R(t)=f_{\rm in}(B_R(t)^c),\qquad
		T_R(t)=\int_{B_R(t)^c}|V_t(z)|^2\dd f_{\rm in}(z).
		\]
		Since $\mathbf1_{\{|X_t|>R\}}\le e^{-aR^\sigma}e^{a|X_t|^\sigma}$, the spatial exponential moment \eqref{eq:spatial-moment} gives, by Markov's inequality,
		\begin{equation}\label{eq:Q}
			Q_R(t)\le M_xe^{-aR^\sigma}\qquad (R>0,\ t\ge0).
		\end{equation}
		Convexity of $u\mapsto|u|^r$ and the alignment equation imply the nonincrease of centered velocity moments:
		\begin{equation}\label{eq:propagate}
			\int|V_t|^r\dd f_{\rm in}\le M_r\qquad(t\ge0).
		\end{equation}
		By H\"older's inequality with exponents $r/2$ and $r/(r-2)$,
		\begin{align}
			T_R(t)
			&\le\left(\int_{B_R(t)^c}|V_t|^r\dd f_{\rm in}\right)^{2/r}
			Q_R(t)^{(r-2)/r}\notag\\
			&\le M_r^{2/r}M_x^{(r-2)/r}
			\exp\!\left(-\frac{a(r-2)}rR^\sigma\right).
			\label{eq:T}
		\end{align}
		Thus, with
		\begin{equation}\label{eq:Db}
			D:=M_r^{2/r}M_x^{(r-2)/r},\qquad b:=\frac{a(r-2)}r>0,
		\end{equation}
		we have $T_R(t)\le De^{-bR^\sigma}$, uniformly in $t$.
		
		\medskip
		\noindent\textbf{Step 2: Localized variance dissipation.}
		The exact variance identity is
		\begin{equation}\label{eq:variance}
			\dot \cV(t)=-\kappa\iint\psi(|X_t-X_t'|)
			|V_t-V_t'|^2\dd f_{\rm in}\dd f_{\rm in}'\le0
			\quad\text{for a.e. }t.
		\end{equation}
		On $B_R(t)\times B_R(t)$ we have $|X_t-X_t'|\le2R$, hence
		$\psi(|X_t-X_t'|)\ge\psi(2R)$. Moreover, since $\int V_t\dd f_{\rm in}=0$,
		\begin{align*}
			\iint_{B_R\times B_R}|V_t-V_t'|^2\dd f_{\rm in}\dd f_{\rm in}'
			&=2(1-Q_R)(\cV-T_R)
			-2\left|\int_{B_R}V_t\dd f_{\rm in}\right|^2\\
			&\ge2(1-Q_R)(\cV-T_R)-2Q_RT_R\\
			&=2\big[(1-Q_R)\cV-T_R\big].
		\end{align*}
		We have therefore proved, for every fixed $R>0$ and almost every $t$,
		\begin{equation}
		\dot \cV(t)\le -2\kappa\psi(2R)
			\big[(1-Q_R(t))\cV(t)-T_R(t)\big].
		\end{equation}
		Choose $R_0\ge1$ so large that $M_xe^{-aR_0^\sigma}\le\frac12$.
		For $R\ge R_0$, \eqref{eq:Q} ensures $1-Q_R(t)\ge\frac12$, and \eqref{eq:T} in \eqref{eq:localized} yields
		\begin{equation}\label{eq:main-ode-radius}
			 \dot \cV(t)\le-\kappa\psi(2R)\cV(t)
				+2\kappa D\psi(2R)e^{-bR^\sigma},\qquad R\ge R_0.
		\end{equation}
		This estimate is valid pointwise in the parameter $R$ at almost every time at which the variance identity holds. Choosing $R=R(t)$ below requires no differentiation of moving spatial sets.
		
		\medskip
		\noindent\textbf{Step 3: Choose a growing radius to balance the two terms.}
		Set
		\begin{equation}\label{eq:thetaR}
			\vartheta:=\frac{\sigma}{\sigma+\beta}\in(0,1),\qquad
			s:=1+t,\qquad
			R(t):=A s^{1/(\sigma+\beta)},
		\end{equation}
		where $A\ge\max\{1,R_0\}$. We have
		$e^{-bR(t)^\sigma}=e^{-Bs^\vartheta}$ with $B:=bA^\sigma$.
		Furthermore, because $s\ge1$,
		\[
		4A^2s^{2/(\sigma+\beta)}\le1+4A^2s^{2/(\sigma+\beta)}
		\le(1+4A^2)s^{2/(\sigma+\beta)},
		\]
		so
		\begin{equation}\label{eq:psi-bounds}
			(1+4A^2)^{-\beta/2}s^{\vartheta-1}
			\le\psi(2R(t))\le
			(2A)^{-\beta}s^{\vartheta-1}.
		\end{equation}
		Using the \emph{lower} bound from \eqref{eq:psi-bounds} for the negative dissipation term and the \emph{upper} bound for the positive error term in \eqref{eq:main-ode-radius}, we obtain
		\begin{equation}\label{eq:scalar}
			\dot \cV(t)\le-k_0s^{\vartheta-1}\cV(t)
				+d_0s^{\vartheta-1}e^{-Bs^\vartheta},
		\end{equation}
		where the constants are explicitly
		\begin{equation}\label{eq:constants}
			k_0:=\kappa(1+4A^2)^{-\beta/2},\quad
			d_0:=2\kappa D(2A)^{-\beta},\quad
			B:=bA^\sigma.
		\end{equation}
		The choice $R(t)\asymp (1+t)^{1/(\sigma+\beta)}$ balances the exponents $tR^{-\beta}$ (cumulative dissipation) and $R^\sigma$ (spatial-tail error).
		
		\medskip
		\noindent\textbf{Step 4: An explicit scalar supersolution.}
		Choose
		\begin{equation}\label{eq:c}
			0<c<\min\left\{B,\frac{k_0}{\vartheta}\right\},
		\end{equation}
		and choose
		\begin{equation}\label{eq:K}
			K\ge\max\left\{e^c\cV(0),\frac{d_0}{k_0-c\vartheta}\right\}.
		\end{equation}
		Define $W(t)=K e^{-c(1+t)^\vartheta}=Ke^{-cs^\vartheta}$. Then $W(0)\ge\cV(0)$ and
		\begin{align*}
			\dot W(t)+k_0s^{\vartheta-1}W(t)
			&=(k_0-c\vartheta)Ks^{\vartheta-1}e^{-cs^\vartheta}\\
			&\ge d_0s^{\vartheta-1}e^{-Bs^\vartheta},
		\end{align*}
		where we used $c<B$ and \eqref{eq:K}. Therefore $W$ is a supersolution of the scalar inequality \eqref{eq:scalar}. The comparison principle gives
		\[
		\cV(t)\le W(t)=K\exp\!\left[-c(1+t)^{\sigma/(\sigma+\beta)}\right],\quad t\ge0,
		\]
		which is \eqref{eq:target}.
	\end{proof}
	
	\begin{proof}[Completion of Theorem~\ref{thm:main}]
		The exponential velocity moment in \eqref{eq:initial-tail} yields
		some moment strictly higher than order two. Apply
		Proposition~\ref{prop:conditional} to the uniform spatial
		tail proved above, obtaining \eqref{eq:main-rate}.
	\end{proof}
	
	\begin{proof}[Proof of Corollary~\ref{cor:W2}]
		The uniform exponential moment implies every polynomial spatial
		moment bound. For each $s$ and large $R$,
		$|X|^s\mathbf1_{|X|>R}\le
		C_s e^{-a_*R^\sigma/2}e^{a_*|X|^\sigma}$,
		which proves uniform integrability. Since $\dot X_t=V_t$,
		\[
		\|X_t-X_s\|_{L^2(f_{\rm in})}\le
		\int_s^t\sqrt{\cV(\tau)}\,\dd\tau.
		\]
		The latter tail integral decays like
		$e^{-c'(1+s)^\vartheta}$, after reducing $c>0$ to absorb
		the polynomial prefactor. Thus $X_t$ converges in $L^2(f_{\rm in})$
		to a random vector $X_\infty$ with the asserted rate.
		The coupling $(X_t,X_\infty)_\#f_{\rm in}$ gives the
		$W_2$ estimate; Fatou's lemma transfers the spatial exponential
		moment to the limit.
	\end{proof}
	
	\section{A countably atomic sharpness theorem}\label{sec:sharpness}
	\begin{theorem}[Optimal spatial-tail-limited exponent for measures]
		\label{thm:atomic}
		Fix $d\ge1$, $\kappa>0$, $0<\beta<1$, and $p>0$.
		There exists a global Lagrangian probability solution of the kinetic
		Cucker--Smale equation, supported on a fixed line in $\mathbb R^d$,
		with bounded initial velocities and unbounded spatial support, such
		that for some $a,c,C>0$,
		\begin{equation}\label{eq:atomic-sharp}
			\sup_{t\ge0}\int e^{a|x-x_c(t)|^p}\,\dd f_t<\infty,
			\qquad
			\mathcal V(t)\ge c\exp\left[-C(1+t)^{p/(p+\beta)}\right],
			\quad t\ge 0.
		\end{equation}
		Hence the exponent $p/(p+\beta)$ in the velocity-variance upper bound
		is optimal for the class of kinetic probability measures having a
		uniform spatial exponential moment of order $p$ and bounded initial
		velocity support. The conclusion concerns atomic measure solutions;
		it does not assert sharpness for smooth strictly positive densities.
	\end{theorem}
\begin{proof}
	It suffices to construct the solution for $d=1$ and embed the
	positions and velocities into one coordinate line for general $d$.
	Set $\alpha=1-\beta$, choose $A>1$ later, and put
	\[
	R_n=A^n,\qquad m_n=\varepsilon e^{-R_n^p}\quad(n\ge1),
	\qquad m_0=1-\sum_{n\ge1}m_n\ge\tfrac12,
	\]
	with $\varepsilon>0$ sufficiently small. Take the symmetric datum
	\begin{equation}\label{eq:atomic-data}
		\mu=m_0\delta_{(0,0)}+
		\sum_{n\ge1}\frac{m_n}{2}
		\big(\delta_{(R_n,1)}+\delta_{(-R_n,-1)}\big).
	\end{equation}
	The central characteristic remains at the origin by symmetry.
	Write $(x_n(t),u_n(t))$ for the positive satellites, with
	$(x_n(0),u_n(0))=(R_n,1)$.
	
	\emph{Global countable dynamics.}
	Introduce $y_n=x_n-R_n$, and regard $(y,u)$ as an unknown in
	$\ell^\infty\times\ell^\infty$, initially $(0,\mathbf1)$.
	The equations are $\dot y_n=u_n$ and
	\begin{equation}\label{eq:atomic-system}
		\dot u_n=-\kappa a_n u_n+
		\frac\kappa2\sum_{j\ge1}m_j
		[\psi(|x_n-x_j|)-\psi(x_n+x_j)]u_j,
	\end{equation}
	where
	\begin{equation}\label{eq:atomic-an}
		a_n=m_0\psi(x_n)+\frac12\sum_{j\ge1}m_j
		[\psi(|x_n-x_j|)+\psi(x_n+x_j)].
	\end{equation}
	Here and below the formula is used while $x_n\ge0$.
	Since $r\mapsto\psi(|r|)$ is bounded and globally Lipschitz,
	the right-hand side, expressed with $x_n=R_n+y_n$ and the
	signed symmetric counterparts, maps
	$\ell^\infty\times\ell^\infty$ to itself and is locally
	Lipschitz on bounded sets, uniformly in $n$: its differences are
	bounded by sums of the form
	$C\sum_jm_j(\|y-\tilde y\|_\infty
	+\|u-\tilde u\|_\infty)$.
	Thus the Banach-space Picard theorem provides a local solution even
	though $(R_n)$ itself is unbounded.
	
	The full symmetric velocity dynamics is a weighted consensus system
	with nonnegative weights. The standard velocity maximum principle
	applied to the bounded sequence gives $|u_n(t)|\le1$; moreover
	$\psi(|x_n-x_j|)\ge\psi(x_n+x_j)$ for $x_n,x_j\ge0$, so
	\eqref{eq:atomic-system} is cooperative on the positive cone.
	Starting from $u_n(0)=1$, its integrating-factor formulation gives
	$u_n(t)\ge e^{-\kappa t}>0$ while all $x_j\ge0$, since
	$a_n\le1$. Consequently $x_n(t)\ge R_n>0$, and the positivity
	bootstrap closes globally. Finally
	$\|u(t)\|_\infty\le1$ and $\|y(t)\|_\infty\le t$, precluding
	finite-time breakdown of the $\ell^\infty$ solution.
	It follows that
	\begin{equation}\label{eq:atomic-monotone}
		0<u_n(t)\le1,\qquad x_n(t)\ge R_n\quad(n\ge1,t\ge0).
	\end{equation}
	
	\emph{Uniform characteristic confinement.}
	Let $\Psi(r)=\int_0^r\psi(|s|)\,\dd s$ be the odd primitive.
	We justify the infinite primitive series appearing in the natural-velocity
	invariant. Since $0<\psi\le1$, we have $|\Psi(s)|\le|s|$, and the
	specific satellite masses satisfy
	\[
	\sum_{j\ge1}m_jR_j
	=\varepsilon\sum_{j\ge1}A^je^{-A^{jp}}<\infty.
	\]
	By \eqref{eq:atomic-monotone},
	$R_j\le x_j(t)\le R_j+t$. In particular, for fixed $n$ and
	$0\le t\le T$,
	\[
	|\Psi(x_n(t)-x_j(t))|
	\le |x_n(t)|+|x_j(t)|
	\le R_n+R_j+2T.
	\]
	The same bound holds for $\Psi(x_n(t)+x_j(t))$, so the
	mass-weighted primitive series (including the reflected satellites)
	converges absolutely and uniformly on each bounded time interval;
	the central-mass contribution is a single finite term. Furthermore,
	\[
	\left|\frac{\dd}{\dd t}\Psi(x_n(t)\mp x_j(t))\right|
	=\psi(|x_n(t)\mp x_j(t)|)\,|u_n(t)\mp u_j(t)|
	\le2.
	\]
	As $\sum_jm_j<\infty$, the series of derivatives converges uniformly
	on bounded time intervals as well. Termwise differentiation is therefore
	legitimate. The characteristic equations now show that the
	natural-velocity quantity
	\[
	u_n(t)+\kappa\int\Psi(x_n(t)-y)\,\dd\rho_t(y)
	\]
	is conserved. For $x,y\ge0$, monotonicity and oddness imply
	\[
	0\le\tfrac12[\Psi(x-y)+\Psi(x+y)]\le\Psi(x).
	\]
	The lower bound follows from $x+y\ge |x-y|$; the upper bound
	follows because the displayed average is nonincreasing as a
	function of $y\ge0$. Using $u_n(t)\ge0$, the central mass, and
	the initial equality $u_n(0)=1$, we obtain
	\begin{equation}\label{eq:atomic-invariant-bound}
		\kappa m_0\Psi(x_n(t))\le1+\kappa\Psi(R_n).
	\end{equation}
	There exist $c_\beta,C_\beta>0$ such that
	$c_\beta r^\alpha\le\Psi(r)\le C_\beta r^\alpha$ for
	$r\ge1$. Since $m_0\ge1/2$ and $R_n\ge1$,
	\eqref{eq:atomic-invariant-bound} implies
	\begin{equation}\label{eq:atomic-conf}
		R_n\le x_n(t)\le C_*R_n,
	\end{equation}
	where $C_*$ depends on $\beta,\kappa$ but not on $n,t,A$ or
	$\varepsilon$ (provided $m_0\ge1/2$).
	
	\emph{Separation and the self-interaction term.}
	Choose $A>4C_*$. By \eqref{eq:atomic-conf}, for $j<n$,
	$|x_n-x_j|\ge 3R_n/4$; for $j>n$,
	$|x_n-x_j|\ge 3R_j/4\ge3R_n/4$.
	Also $x_n+x_j\ge R_n$. Thus, for $j\ne n$,
	\[
	\psi(|x_n-x_j|)+\psi(x_n+x_j)\le C R_n^{-\beta}.
	\]
	The $j=n$ summand in \eqref{eq:atomic-an} has
	$\psi(0)=1$ and contributes at most $m_n/2$ plus
	$C m_nR_n^{-\beta}$. Since
	$m_n\le\varepsilon e^{-R_n^p}\le C_pR_n^{-\beta}$,
	this term is also absorbed. Hence
	\begin{equation}\label{eq:atomic-damp}
		0\le a_n(t)\le C_0R_n^{-\beta}
		\qquad(n\ge1,t\ge0),
	\end{equation}
	with $C_0$ independent of $n,t$. The nonnegative off-diagonal
	term in \eqref{eq:atomic-system} therefore gives
	\begin{equation}\label{eq:atomic-u-lower}
		u_n(t)\ge\exp(-\kappa C_0tR_n^{-\beta}).
	\end{equation}
	
	\emph{Uniform tails and the lower rate at every large time.}
	The zero mean velocity implies
	$\mathcal V(t)=\sum_{n\ge1}m_nu_n(t)^2$, and hence
	\begin{equation}\label{eq:atomic-individual-lower}
		\mathcal V(t)\ge\varepsilon\exp
		\big[-R_n^p-2\kappa C_0tR_n^{-\beta}\big]
		\quad\text{for each }n\ge1.
	\end{equation}
	For $t\ge A^{p+\beta}$ choose the integer $n\ge1$ satisfying
	\[
	A^n\le t^{1/(p+\beta)}< A^{n+1}.
	\]
	Then $A^{-1}t^{1/(p+\beta)}<R_n\le t^{1/(p+\beta)}$, so
	both terms in the exponent of
	\eqref{eq:atomic-individual-lower} are bounded by
	$C(1+t)^{p/(p+\beta)}$. Enlarging $C$ and decreasing the
	prefactor covers bounded times as well (the variance is positive
	there). This proves the second assertion in \eqref{eq:atomic-sharp}.
	Finally, for $0<a<C_*^{-p}$, \eqref{eq:atomic-conf} yields
	\[
	\sup_{t\ge0}\int e^{a|x|^p}\,\dd f_t
	\le m_0+\sum_{n\ge1}m_ne^{aC_*^pR_n^p}<\infty.
	\]
	This proves the uniform exponential spatial moment and the theorem.
\end{proof}

	\begin{remark}[Infinite-particle background]
		Ha, Wang, and Xue \cite{HaWangXueInfinite2026} study countably
		infinite Cucker--Smale systems using infinite-norm inequalities and
		spectral criteria for clustering. The construction here instead uses
		a summably weighted symmetric atomic measure and the bounded
		\(\ell^\infty\) displacement--velocity system; its conclusion is a
		lower bound for the velocity variance of a kinetic probability
		solution, rather than a spectral characterization of clusters.
	\end{remark}
	
	\section{Discussion and open problems}
	\label{sec:discussion}
	The finite-moment theorem gives uniform spatial second moments,
	strong confinement along characteristics, integrable velocity
	fluctuations and a polynomial variance-decay estimate. Exponential
	initial position and velocity moments improve this to uniform
	exponential spatial tails, stretched-exponential velocity alignment,
	and convergence of the centered spatial law in $W_2$.
	The countably atomic construction proves the spatial-tail-limited
	rate exponent optimal in the measure-valued class, with bounded
	initial velocities. We do not infer sharpness in the smooth-density
	class or in the velocity-limited regime from that construction.
	
    Several questions remain open and suggest directions for further investigation. First, it would be interesting to determine whether the finite-moment confinement theorem extends to the critical velocity-moment exponent, i.e. $m=2/(1-\beta)$. A second direction concerns the optimality of the quantitative alignment estimates. Although our countably atomic construction establishes sharpness in the spatial-tail-limited regime, it remains unclear whether the same obstruction can be realized by smooth, strictly positive phase-space densities, or whether the predicted decay rate is optimal when the initial velocity tails determine the convergence exponent. Finally, an important challenge is to extend the confinement and tail-propagation results to critical communication decay, particularly the borderline case where the interaction weight decays inversely with distance i.e. $\beta=1$. The borderline communication case $\beta=1$ has already been treated conditionally by Ha and Wang, who establish weak flocking with uniformly bounded centered spatial second moments and algebraic velocity alignment under an initial exponential spatial moment, bounded velocity support, and an explicit lower bound on the coupling strength. The question here is therefore not whether weak flocking can occur at $\beta=1$, but whether the support-independent characteristic confinement and strong spatial-tail propagation obtained above can be extended to that critical regime with unbounded velocity support and without such a coupling restriction, as well as to other borderline long-range communication kernels.
	
	\section{Declaration of Generative AI Use}
	
	During the preparation of this manuscript, the author used ChatGPT (OpenAI) as a research-assistance tool for exploring mathematical arguments, examining proof strategies, checking calculations, reviewing relevant literature, and improving the organization and clarity of the text. AI-generated suggestions and derivations were treated as material requiring mathematical verification, rather than as independent evidence of correctness. The author retains full responsibility for the validity of all mathematical statements, the accuracy of the references, and the final content of the manuscript.

\end{document}